\documentclass[a4paper,12pt]{amsart}
\usepackage{amsfonts}
\usepackage{bbm}
\usepackage{hyperref}
\usepackage{mathrsfs}

\newtheorem{theorem}{Theorem}[section]
\newtheorem{lemma}[theorem]{Lemma}

\theoremstyle{definition}

\theoremstyle{remark}

\numberwithin{equation}{section}

\newcommand{\hsp}{\hspace*{\parindent}}

\begin{document}

\title[]{\large A note on the Gauss-Bonnet-Chern theorem for general connection}

\author{Haoyan Zhao}

\address{ Academy of math and systems sciences, Chinese academy of science,
Beijing, 100190
\\\hsp
Mathematical Sciences center, Zhejiang University, Hangzhou, 310028, People's Republic of China}
\email{zhaohaoran81@163.com}

\subjclass[2000]{58J20}

\date{}

\keywords{heat equation, Hodge-Laplacian}

\begin{abstract}In this paper, we prove a local index theorem for the DeRham Hodge-Laplacian which is defined by the connection compatible with metric.
This connection need not be the Levi-Civita connection. When the connection is Levi-Civita connection, this is the classical
local Gauss-Bonnet-Chern theorem.
\end{abstract}

\maketitle

\section{Introduction}

The Gauss-Bonnet-Chern theorem is an index theorem about relationship between topology and geometry on a compact manifold. It has been proved by Allendoerfer-Weil \cite{al} in 1940s.
Later, Chern \cite{chen} has given a proof by intrinsic computation. The refined local Gauss-Bonnet-Chern theorem was proved by Patodi \cite{patodi} in 1971, which was conjectured by Mckean-Singer \cite{mckean}.
In the above theorem, the DeRham operator and Hodge-Laplacian are defined by Levi-Civita connection. Recently, Beneventano-Gilkey-Kirsten-Santangelo \cite{ben} have studied the Gauss-Bonnet theorem for general connection and corresponding heat trace's asymptotic expansion.
Bell \cite{bell} has given Gauss-Bonnet theorem for vector bundle whose rank is equal to the dimension of underlying
manifold. There have been also some works about Gauss-Bonnet-Chern
theorem's generalization in Finsler geometry (see Bao-Chern-Shen \cite{bao}, Lackey \cite{lackey}, Zhao \cite{zhao} and so on).
\\\hsp Now we state the main theorem in this article.
\begin{theorem}\label{theorem} Let M be a compact Riemannian manifold of even dimension $d$ $(d=2l)$ with metric g, which has a metric compatible connection D. Let
$\varepsilon(e)$ denote exterior multiplication by differential form $e$ and
$\iota(e)$ denote interior multiplication by $e$. Then we may define
DeRham Dirac operator, for exterior form bundle section f:
\begin{equation*}\mathcal {D}f=(\varepsilon(e^{i})-\iota(e^{i}))D_{e_{i}}f\end{equation*}
Where $\{e_{i}\}$is tangent vector frame, $\{e^{i}\}$ is its dual frame.
Let h(t,x,y) be the heat kernel (fundamental solution) for the following DeRham Hodge-laplacian equation:
\begin{equation*}\frac{\partial}{\partial t}f=-\frac{1}{2}\mathcal {D}^{2}f\end{equation*}

Then we have
\begin{equation*}\lim_{t\rightarrow0}Str\{h(t,x,x)\}dm= \frac{1}{(2\pi)^{l}}Pf(-R),   \hsp (1.1)
\end{equation*}

where $Str$ denote supertrace, $dm$ is the volume element, $R$ is the Riemannian curvature corresponding to $D$, $Pf$ is Pfaffian.

\end{theorem}

\section*{Acknowledgement}
 The author is  grateful to Chinese academy
science and its academy of math and systems science, math institute, math and interdisciplinary science center, morningside center; math
science center and math department in Zhejiang university; math science center and math department in
Tsinghua university for their supports and help. The author also thanks Professor Youde Wang , Professor Kefeng Liu and Professor Xiangyu Zhou for
their supports and help . The author also thanks all teachers and friends.

\section{geometric preliminares}

Let $\widehat{D}$ denote Levi-Civita connection, $\{e_{i}\}$ be
orthonormal tangent frame.
$\widetilde{D}(e_{i})={\hat{\Gamma}}_{li}^{s}e^{l}\otimes e_{s}$, $(D)(e_{i})=\Gamma_{li}^{s}e^{l}\otimes e_{s}$.
Let $$\wedge T^{*}_{x}M=\wedge ^{+}T^{*}_{x}M\oplus\wedge ^{-}T^{*}_{x}M,$$ for $x\in M$, where $\wedge ^{+}T^{*}_{x}M$ consists of even degree forms, $\wedge ^{-}T^{*}_{x}M$ consists of odd degree forms.
For $a\in End(\wedge T^{*}_{x}M)$, define $$Str(a)=trace(a_{\wedge^{+}})-trace(a_{\wedge^{-}}).$$

The fundamental solution's asymptotic expansion in $(x,x)$ is determined by local condition around $x$ (see [4] chapter 2). So we may assume $M$ is spin manifold,
whose spinor bundle is denoted by $\mathcal {S}$, dual bundle is denoted by $\mathcal {S}^{*}$.
For connection $\widehat{D}$ and $D$ respectively, $\mathcal {S}$ and $\mathcal{S}^{*}$ have lifted connections, which are  denoted by
$\widehat{D}^{\mathcal {S}}$ and $\widehat{D}^{\mathcal {S}^{*}}$, $D^{\mathcal {S}}$ and $D^{\mathcal {S}^{*}}$.
There  is  a linear isomorphism between Clifford bundle $Cl(T^{*}M)$ and exterior form bundle $\Omega^{*}(M)$, which is not an
algebraic isomorphism:
$$Cl(T^{*}M)\cong\Omega^{*}(M).\hsp(1.2)$$
The isomorphism (1.2) is denoted by by $c$ , the inverse of $c$ by $\sigma$   .
\\\hsp We still denote complexified Clifford algebra bundle
by $Cl(T^{*}M)$. $Cl(T^{*}M)$ has an action on $\mathcal {S}$. The
lifted connection on $\mathcal {S}$ is compatible with this Clifford
action. There is also an isomorphism between
the complexified Clifford algebra bundle $Cl(T^{*}M)$ and the endmorphism bundle $\mathcal {S}\otimes\mathcal {S}^{*}$ of $\mathcal {S}$: $$Cl(T^{*}M)\cong\mathcal {S}\otimes\mathcal {S}^{*}.\hsp(1.3)$$ This isomorphism is an algebraic
isomorphism.
For a orthonormal base $\{e^{i}\}\in T^{*}_{x}M$, we define a
chirality
element: $$\Gamma=i^{l}c(e^{1})c(e^{2})...c(e^{d}).$$
By computations , one could get $\Gamma^{2}=1$. Let $$\mathcal {S}^{+}=\{a|\Gamma a=a, a\in
\mathcal {S}_{x}\},$$ $$\mathcal {S}^{-}=\{a|\Gamma a=-a, a\in
\mathcal {S}_{x}\}.$$

By isomorphism $(1.2)$ and $(1.3)$, there is:
$$\Omega^{*}(M)\cong\mathcal {S}\otimes\mathcal {S}^{*}.\hsp(1.4)$$
As exterior form bundle, $\Omega^{*}(M)$ has a connection which is
from  $T^{*}M$'s connection, a graded structure in terms of even and
odd degree. $Cl(T^{*}M)$ has an action on it. As tensor product of $\mathcal {S}$ and $\mathcal
{S}^{*}$, $\mathcal {S}\otimes\mathcal {S}^{*}$ has  an
connection from $\mathcal {S}$, a graded structure $$((\mathcal {S}^{+}\otimes\mathcal
{S}^{*+})\oplus (\mathcal{S}^{-}\otimes\mathcal {S}^{*-}))\bigoplus((\mathcal {S}^{+}\otimes\mathcal
{S}^{*-})\oplus (\mathcal {S}^{-}\otimes\mathcal {S}^{*+})).$$
$Cl(T^{*}M)$ also has an action on it. Under the above isomorphism
(1.4), these two connections, graded structure, action are identical (see Berline-Getzler-Vergne \cite{berline} Chapter 3 and 4).
So from now we will always consider exterior form bundle $\Omega^{*}(M)$ as twisted Clifford module bundle $\mathcal {S}\otimes\mathcal
{S}^{*}$.
\begin{lemma}\label{2}

Let $T$ denote Berezin integral (see Berline-Getzler-Vergne \cite{berline}). For $a\in End(\mathcal {S}_{x})\cong Cl(T^{*}_{x}M)$, $b\in End(\mathcal {S}^{*}_{x})\cong Cl(T^{*}_{x}M)$, $$End(\mathcal{S}_{x})\otimes End(\mathcal {S}^{*}_{x})\cong End(\wedge^{*}(T^{*}_{x}M\otimes_{R}
C)),$$
$$Str(a)=tr(\Gamma a)=(-2i)^{l}T\circ\sigma(a),$$
$$Str(b)=tr(\Gamma^{*}b)=(2i)^{l}T\circ\sigma(b),$$
$$Str(a\otimes b)=tr(\Gamma a)tr(\Gamma^{*}b).$$
\end{lemma}
For the Dirac operator
$ \mathcal {D}^{\mathcal {S}}$ associated with $D^{\mathcal {S}}$:$$ \mathcal {D}^{\mathcal {S}}=c(e^{i})D^{\mathcal {S}}_{e_{i}},$$ there exist unique 1-form $a=a_{i}e^{i}$, 3-form$B=B_{il}^{s}e^{i}\wedge e^{l}\wedge e^{s}$, such that
$$\mathcal {D}^{\mathcal {S}}=\widehat{\mathcal {D}}^{\mathcal {S}}+c(a)+c(B).$$
Define
$$D^{\mathcal {S},B}_{e_{i}}=\widehat{D}^{\mathcal {S}}_{e_{i}}+B_{il}^{s}c(e^{l})c(e^{s}),$$
then $\mathcal {D}^{\mathcal {S},B}=\widehat{\mathcal {D}}^{\mathcal {S}}+c(B)$.
\\\hsp Let $W$ be a complex bundle equipped with connection $D^{W}$, curvature be
$F$. Define the connections
$\widehat{D}^{\mathcal {S}\otimes W}, D^{\mathcal {S}\otimes W} , D^{\mathcal {S}\otimes W,B},D^{\mathcal {S}\otimes W,3B}$, whose corresponding Dirac operator are $\widehat{\mathcal
{D}}^{\mathcal {S}\otimes W},\mathcal {D}^{\mathcal {S}\otimes W},\mathcal {D}^{\mathcal {S}\otimes W,B},\mathcal {D}^{\mathcal {S}\otimes W,3B}$:
$$ \widehat{D}^{\mathcal {S}\otimes W}=\widehat{D}^{\mathcal {S}}\otimes1+1\otimes D^{W},$$
$$ D^{\mathcal {S}\otimes W}=D^{\mathcal {S}}\otimes1+1\otimes D^{W},$$
$$ D^{\mathcal {S}\otimes W,B}=D^{\mathcal {S},B}\otimes1+1\otimes D^{W},$$
$$ D^{\mathcal {S}\otimes W,3B}=D^{\mathcal {S},3B}\otimes1+1\otimes D^{W}.$$

A useful formula on the square of $\mathcal {D}$ is from Bismut \cite{bis1}.

\begin{lemma}\label{3}(Bismut \cite{bis1})$$(\mathcal {D}^{\mathcal {S}\otimes W,B})^{2} =-\triangle^{\mathcal {S}\otimes W,3B} +\frac{s}{4}
+c(F)+c(dB)-2|B|^{2} .$$
\begin{equation*}
\begin{split} (\mathcal{D}^{\mathcal {S}\otimes W}) ^{2} &=  (\mathcal {D}^{\mathcal {S}\otimes W,B})^{2}-2(a,e^{i}
)\widehat{D}^{\mathcal {S}\otimes W}_{e_{i}}+c(\widehat{D}a)-2c(\iota(a)B)-|a|^{2} \\&=  -\triangle^{\mathcal {S}\otimes W,3B}- 2(a,e^{i})D^{\mathcal {S}\otimes W,3B}_{e_{i}}
+C,
\end{split}
\end{equation*}
where $C=
\frac{s}{4}+c(F)+c(dB)-2|B|^{2}+c(\widehat{D}a)-|a|^{2}$, $s$ is scalar curvature.\end{lemma}
Define the connections $\widehat{D}^{\mathcal {S}\otimes \mathcal {S}^{*}}, D^{\mathcal {S}\otimes \mathcal {S}^{*}}, D^{\mathcal {S}\otimes \mathcal {S}^{*},B}, D^{\mathcal {S}\otimes \mathcal {S}^{*},3B}$ on $\mathcal {S}\otimes\mathcal
{S}^{*}$, and corresponding Dirac operators are noted by $\widehat{\mathcal
{D}}^{\mathcal {S}\otimes \mathcal {S}^{*}},\mathcal {D}^{\mathcal {S}\otimes \mathcal {S}^{*}},
\mathcal {D}^{\mathcal {S}\otimes \mathcal {S}^{*},B},\\\mathcal {D}^{\mathcal {S}\otimes \mathcal {S}^{*},3B}$:
$$ \widehat{D}^{\mathcal {S}\otimes \mathcal {S}^{*}}=\widehat{D}^{\mathcal {S}}\otimes1+1\otimes D^{\mathcal {S}^{*}},$$
$$ D^{\mathcal {S}\otimes \mathcal {S}^{*}}=D^{\mathcal {S}}\otimes1+1\otimes D^{\mathcal {S}^{*}},$$
$$ D^{\mathcal {S}\otimes \mathcal {S}^{*},B}=D^{\mathcal {S},B}\otimes1+1\otimes D^{\mathcal {S}^{*}},$$
$$ D^{\mathcal {S}\otimes \mathcal {S}^{*},3B}=D^{\mathcal {S},3B}\otimes1+1\otimes D^{\mathcal {S}^{*}}.$$
When consider $ \Omega^{*}(M)$ as $\mathcal {S}\otimes\mathcal
{S}^{*}$, $ D^{\mathcal {S}\otimes \mathcal {S}^{*}}$ and $\mathcal {D}^{\mathcal {S}\otimes \mathcal {S}^{*}}$ are respectively $D$ and $\mathcal {D}$ defined
in theorem \ref{theorem}. So we can get the expression of $\mathcal {D}^{2}$ by lemma
\ref{3}. This is the key step.
\section{the proof of the main theorem}
In the following proof, we use Feynman-Kac formula and the generalized Wiener functional which
are
included in stochastic analysis. In the course of studying Malliavin theory, the generalized Wiener functional
and its applications were introduced and studied by Malliavin, Kusuoka-Stroock \cite{stroock1}, Watanabe \cite{watanabe1} \cite{watanabe2}, Ikeda-Watanabe \cite{ikeda-watanabe} and so on.
In this paper we adopt the definition and processing mode as in Watanabe \cite{watanabe2}. Watanabe \cite{watanabe2} proved the local
Gauss-Bonnet-Chern theorem and signature theorem by this
method. About more details on generalized Wiener functional and its
applications, the readers could refer to Ikeda-Watanabe \cite{ikeda-watanabe}. The
probabilistic proof on index theorem was provided firstly by Bismut \cite{bis2}.
Base on probabilistic method , Bismut \cite{bis1} proved a local index theorem
on non K$\ddot{a}$hler manifold. There was still a stochastic proof for the local Gauss-Bonnet-Chern theorem in Elton Hsu \cite{hsu}. If not using stochastic method ,
the main
theorem in this article should be also able to  be proved by Getzler's rescaling method as in Berline-Getzler-Vergne \cite{berline} and
Getzler \cite{getzler}. The most notations in our computations are the same as
Watanabe \cite{watanabe2}.
\\\hsp Furthermore, we assume $M$ be $R^{d}$, with metric $g$, which
is equal to the standard Euclidean metric outside of some sufficient
big ball, the natural coordinate on $R^{d} $ be identical to normal coordinate around the original point
.

Let$\{e_{i}\}$ be natural tangent frame, $\{f_{i}\},\{g_{i}\}$  be orthonormal frames respectively by parallel translations along
the radial lines from original point under connections $D^{3B},
D$ on $M$, $ \{e^{i}\}, \{f^{i}\},\{g^{i}\}$ are respectively their dual frames. Let connection $D^{3B}=({\Gamma}^{3B})_{il}^{s}e^{i}\otimes f_{s}\otimes f^{l}$,
under frame $\{f_{i}\}$ ;  connection
$D={\Gamma}^{s}_{il}e^{i}\otimes g_{s}\otimes g^{l}$ under frame
$\{g_{i}\}$. By virtue of the above frames, we could get the
trivialization of $\mathcal {S}\otimes\mathcal {S}^{*}$: $R^{d}\times (S \otimes
S^{*})$, where $S$ is the spinor space of Euclidean space $R^{d}$.

Let   $C_{i}=\frac{1}{4}({\Gamma}^{3B}) _{il}^{s} c(f^{l}
)c(f^{s})\otimes1+1\otimes(-\frac{1}{4}{\Gamma}_{il}^{s}) c^{*}(g^{s}) c^{*}(g^{l})$, $b^{i}=-g^{ls}{\hat{\Gamma}}
_{ls}^{i}+2(a,e^{i})$, then according to lemma \ref{3}, the
Hodge-Laplacian heat equation expression in the natural coordinate is:
\[\left\{
        \begin{array}{ll}
            \displaystyle\frac{\partial f}{\partial
             t}=\frac{1}{2}g^{ij}(\frac{\partial}{\partial x^{i}}+C_{i})(\frac{\partial}{\partial x^{j}}+C_{j})f
             +\frac{1}{2}b^{i}(\frac{\partial}{\partial x^{i}}+C_{i})f-\frac{1}{2}Cf,~~(t,x)\in(0,\infty)\times R^{d}\\
\displaystyle f(0,x)=\varphi(x),~~ x\in R^{d}
        \end{array}
\right.
\]

Let  smooth real symmetric positive matrix $\sigma^{i}_{k}$  make $$\sum_{k}\sigma^{i}_{k}
\sigma^{l}_{k}=g^{il},$$
then consider the following stochastic differential equation valued in   $R^{d}\times (Cl(R^{d})\otimes Cl(R^{d}))\times(Cl(R^{d})\otimes Cl(R^{d}))$,

\[\left\{\begin{array}{ll}\displaystyle dX^{i}(t)=\sigma^{i}_{k}(X(t))dw^{k}(t)+ \frac{1}{2}b^{i}(X(t))dt, \\
        \displaystyle de(t)=e(t)C_{i}(X(t))\circ dX^{i}(t),\\
        \displaystyle dM(t)=-\frac{1}{2}M(t)e(t)C(X(t))e^{-1}(t)dt,\\
         \displaystyle e(0)=1, M(0)=1, X(0)=x, \end{array}
\right.
\]
$e(t)$ are inverse almost everywhere (see Stroock \cite{stroock}), so $e^{-1}(t)$ are well
defined.

By It$\bar{o}$ formula and properties of generalized Wiener functional ,
$$f(t,x)=E[M(t)e(t)\varphi(X_{x}(t))],$$
$$h(t,x,y)=E[M(t)e(t)\delta_{y}(X(t))],$$
where $\delta_{y}$ is a generalized function: the Dirac delta function associated with $y$.
 It is difficult that compute directly asymptotic expansion for $t$ according to the above formula.  As in Bismut \cite{bis2}, Watanabe \cite{watanabe2}, we consider the stochastic
differential equations with parameter $\varepsilon$:

 \[\left\{
         \begin{array}{ll}
             \displaystyle dX^{i}(t)=\varepsilon\sigma^{i}_{k}(X(t))dw^{k}(t)-\frac{\varepsilon^{2}}{2}b^{i}(X(t))dt, \\
 \displaystyle de(t)=e(t)C_{i}(X(t))\circ dX^{i}(t),
 \hsp i,j,k=1,2,...d. \displaystyle\\(X(0),e(0) =(0,1),
 \end{array} \right.
 \]
 \[\left\{
         \begin{array}{ll}
             \displaystyle dM(t)=-\frac{\varepsilon^{2}}{2}M(t)e(t)C(X(t))e^{-1}(t)dt,\\
 \displaystyle M(0)=1,
   \end{array}
 \right.
 \]
 let us denote the solution  by
 $r^{\varepsilon}(t)=(X^{\varepsilon}(t),e^{\varepsilon}(t),M^{\varepsilon}(t))$, then
$$h(\varepsilon^{2},0,0)=E[M^{\varepsilon}(1)e^{\varepsilon}(1)\delta_{0}(X^{\varepsilon}(1))].$$
\begin{lemma}\label{4}(see Berline-Getzler-Vergne \cite{berline})$$\Gamma_{i}(x)=-\frac{1}{2}\sum_{j}R(\partial_{i},\partial_{j})(0)x^{j}+O(|x|^{2}).$$
$$\Gamma^{3B}_{i}(x)=-\frac{1}{2}\sum_{j}R^{3B}(\partial_{i},\partial_{j})(0)x^{j}+O(|x|^{2}).$$

\end{lemma}

\begin{lemma}\label{1}(see Watanabe \cite{watanabe1} \cite{watanabe2}) Let $D^{\infty}$ be the
 space consists of $R^{d}$ valued Wiener functionals whose any order Malliavin
 derivatives
 are $L_{p}$ integrable, for all
 $p>1$, $\widetilde{D}^{-\infty}$ be its dual space.

 $$X^{\varepsilon}(1)=\varepsilon w(1)+O(\varepsilon^{2})$$~~~~~~~~~in $D^{\infty}$.

 $$\delta_{0}(X^{\varepsilon}(1))=\varepsilon^{-d}\delta_{0}(w(1))+O(\varepsilon^{-d+1})$$
in $\widetilde{D}^{-\infty}$.

 $$E[\delta_{0}(w(1))\cdot
\Phi(w)]=(2\pi)^{-l}E[\Phi(w)|w(1)=0],\hsp
\Phi\in\widetilde{D}^{\infty}.$$

\end{lemma}
Let
  \begin{eqnarray*}\theta^{\varepsilon}(t)&=&\int^{t}_{0}
 C_{i}(X^{\varepsilon}(t))\circ d(X^{\varepsilon})^{i}(s)\\
 &=&\varepsilon^{2}(C^{1}_{ij}(t)c(f^{i})c(f^{j})\otimes1+1\otimes
 C^{2}_{ij}(t)c^{*}(g^{i})c^{*}(g^{j}))+O(\varepsilon^{3}),\end{eqnarray*}
  in which
 $$C^{1}_{ij}(t)=\frac{1}{8}R^{3B}_{mkij}(0)\int^{t}_{0}w^{k}(s)\circ d w^{m}(s),$$
 $$C^{2}_{ij}(t)=\frac{1}{8}R_{mkij}(0)\int_{0}^{t}w^{k}(s)\circ dw^{m}(s),$$

 $ R_{mkij}=(R(e_{m},e_{k})g_{j},g_{i}), R^{3B}_{mkij}=(R^{3B}(e_{m},e_{k})f_{j},f_{i}).$
\\let$$B^{\varepsilon}(t)=\theta^{\varepsilon}(t)-\int_{0}^{t}\frac{\varepsilon^{2}}{2}C(X^{\varepsilon}(s))ds.$$
 \begin{eqnarray*}M^{\varepsilon}(1)e^{\varepsilon}(1)&=&1+\int^{1}_{0}M^{\varepsilon}(s)e^{\varepsilon}(s)\circ
 d\theta^{\varepsilon}(s)+\int^{1}_{0}M^{\varepsilon}(s)e^{\varepsilon}(s)(-\frac{\varepsilon^{2}}{2}C(X^{\varepsilon}(s)))ds
 \\&=&1+B^{\varepsilon}(1)+\int^{1}_{0}\int^{t_{1}}_{0}M^{\varepsilon}(t_{2})e^{\varepsilon}(t_{2})\circ dB^{\varepsilon}(t_{2})
 \circ dB^{\varepsilon}(t_{1})
 \\&=&1+B^{\varepsilon}(1)+\int^{1}_{0}B^{\varepsilon}(t_{1})\circ dB^{\varepsilon}(t_{1})\\&+&\int^{1}_{0}
 \int^{t_{1}}_{0}\int^{t_{2}}_{0}M^{\varepsilon}(t_{3})e^{\varepsilon}(t_{3})\circ dB^{\varepsilon}(t_{3})\circ dB^{\varepsilon}(t_{2})
 \circ dB^{\varepsilon}(t_{1})
 \\
 &=&1+A_{1}+A_{2}+...+A_{l}+O(\varepsilon^{2l+2})
 \end{eqnarray*} in $D^{\infty}( Cl(R^{d})\otimes End(R^{s}))$,\\in which £¬
 \begin{eqnarray*}A_{m}&=&\int^{1}_{0}\int^{t_{1}}_{0}\int^{t_{2}}_{0}...\int^{t_{m-1}}_{0}\circ dB^{\varepsilon}
 (t_{m})\circ dB^{\varepsilon}(t_{m-1})\circ...\circ dB^{\varepsilon}(t_{1})
 \\&=&\varepsilon^{2m}\int^{1}_{0}\int^{t_{1}}_{0}...\int^{t_{m-1}}_{0}\circ d \widetilde{C}(t_{m})\circ d \widetilde{C}(t_{m-1})\circ...\circ
  d \widetilde{C}(t_{1})+O(\varepsilon^{2m+1})\end{eqnarray*}
  in $D^{\infty}( Cl(R^{d})\otimes End(R^{s}))$,  $$\widetilde{C}(t)=C^{1}_{ij}(t)c(f^{i})c(f^{j})\otimes1+1\otimes
 C^{2}_{ij}(t)c^{*}(g^{i})c^{*}(g^{j})-\int_{0}^{t}\frac{1}{2}C(0)dt.$$

Note lemma \ref{2}, when $m<l,$
$$Str(A_{m})=0,$$

\hsp$m=l,$
$$Str(A_{m})=Str(A_{l}),$$
\hsp$m>l$,
$$Str(A_{m})=O(\varepsilon^{2l+3})=O(\varepsilon^{d+3}).$$
\begin{eqnarray*}
Str(A_{l})&=&\frac{\varepsilon^{2l}}{l!}Str\{(-\frac{1}{2}C(0))^{l}\}+O(\varepsilon^{2l+1})
\\&=&\frac{\varepsilon^{2l}}{l!}Str\{(-\frac{1}{2}c(F)(0))^{l}\}+O(\varepsilon^{2l+1})
\\&=&\frac{\varepsilon^{2l}}{l!}Str\{(-\frac{1}{4}(\frac{1}{4}R^{}_{ijnm}(0)c(g^{i})c(g^{j})c^{*}(g^{m})c^{*}(g^{n})))^{l}\}+O(\varepsilon^{2l+1})
\\&=&\varepsilon^{2l}(-2i)^{2l}(-1)^{l} Pf(-\frac{1}{4}R^{}(0))+O(\varepsilon^{2l+1})
\\&=&\varepsilon^{2l}Pf(-R^{}(0))+O(\varepsilon^{2l+1})
\end{eqnarray*}

Note lemma \ref{1}, when $\varepsilon\rightarrow0$£¬\begin{eqnarray*} &&Str[h(\varepsilon^{2},0,0)]e^{1}\wedge e^{2}\wedge...\wedge e^{d}
 =\frac{1}{(2\pi)^{l}}Pf(-R^{}(0))+O(\varepsilon).\end{eqnarray*}\\
 Therefore, we get

 \hsp  $$\displaystyle\lim_{t\rightarrow0}Str
 [h(t,0,0)]e^{1}\wedge e^{2}\wedge...\wedge e^{d}=\frac{1}{(2\pi)^{l}}Pf(-R^{}(0)).$$


\end{document}